\documentclass[12pt]{article}
\usepackage{amsmath}
\usepackage{amsfonts}
\usepackage{amssymb}
\usepackage{theorem}
\usepackage{xcolor}
\usepackage{enumerate}

\title{Subgroups of $SF(\omega )$ and the relation of almost containedness}
\author{B.Majcher-Iwanow 
} 
\date{ } 

\newtheorem{THM}{Theorem}[section]
\newtheorem{LEM}[THM]{Lemma}

\newtheorem{PRO}[THM]{Proposition}

\def\N2N{\omega^\omega}
\def\22O{2^\omega}

\def\PART{(\omega)}
\def\FPART{(\omega)^{<\omega}}

\def\LATTICE{(\PART,\leq)}

\begin{document}
\maketitle
\topskip 20pt

\begin{quote}
{\bf Abstract.} 
The relations of almost containedness and 
othogonality in the lattice of groups of 
finitary permutations are studied in the paper.
We define six cardinal numbers naturally 
corresponding to these relations by 
the standard scheme 
of $P(\omega)$. 
We obtain some consistency results 
concerning these numbers and some versions 
of the Ramsey theorem.

\bigskip


\end{quote}

\bigskip

\setcounter{section}{0}
\section{Preliminaries}

\subsection{Introduction}

The paper is motivated by investigations  
of various versions of van Douwen's diagram,  
i.e. the set of relations between  six  cardinals  referring
to  simple  properties of almost disjointness and almost
containedness, for example see 
\cite{brendle1}, \cite{brendle2}, \cite{brendle3}, 
\cite{c}, \cite{vd}, \cite{minami}, \cite{vaughan}.  
The following theorem proved by
P. Matet in \cite{mat} became one of the motivating results 
in this direction:
\begin{quote}
Let $(\omega)^{\omega}$ be the set of all partitions of $\omega$
having infinitely many classes. Let $\le$ be the order on
$(\omega)^{\omega}$ defined by: $E_{1} \le E_{2}$ if $E_{2}$ is
finer than $E_{1}$. Then assuming the continuum hypothesis there
is a filter $F \subset (\omega)^{\omega}$ such that for every
$({\bf \Sigma}^{1}_{1} \cup  {\bf \Pi}^{1}_{1})$-coloring
$\delta : (\omega)^{\omega} \rightarrow 2$ there is a partition
$E \in F$ such that $\delta$ is constant on the set of all
infinite partitions coarser than $E$.
\end{quote}

The statement is a variant (and a consequence) of the
dualized version of Ramsey's theorem proved by T. Carlson 
and S. Simpson in \cite{b}. 
The argument of P. Matet uses the observation 
that the tower cardinal (we denote it by ${\mathfrak t}_{d}$)
for the ordering of infinite partitions is uncountable. 
Here ${\mathfrak t}_{d}$ is defined by the same scheme as 
${\mathfrak t}$ for the lattice $P(\omega)$ of all 
subsets of $\omega$ (see \cite{vd}).
Moreover it is proved in \cite{mat} that 
the tower cardinal for partitions is $\omega_{1}$  
in {\bf ZFC} and it is proved in \cite{c} that 
the size of a maximal almost orthogonal family 
of partitions must be $2^{\omega}$.
\parskip0pt

The lattice of partitions under the order reversing 
$\le$ was studied in \cite{ma3} and \cite{brendle2}. 
It is shown there that the corresponding 
cardinal invariants look differently:
each of them, except ${\mathfrak a}$, is equal to the analogous
cardinal  in the lattice $P(\omega)$.

Note that this lattice can be also defined 
to be the lattice of 1-closed subgroups of $Sym(\omega)$ 
(i.e. the automorphism groups of structures 
with  unary predicates only). 
Indeed, the corresponding isomorphism of these 
lattices maps a partition $E$ to the group $G_{E}$
of all permutations preserving the $E$-classes.
On the other hand, for any $E$ the group $G_{E}$
is uniquely determined by the subgroup 
$G_{E} \cap SF(\omega)$ of the group $SF(\omega)$ 
of all finitary permutations of $\omega$.
The embedding obtained 
$E \rightarrow G_{E} \cap SF(\omega)$ maps
almost trivial equivalence relations into 
the ideal $IF$ of all finite subgroups of $SF(\omega)$. 

This motivates further questions. 
For example it is interesting to find 
a variant of the result of P. Matet in the lattice 
of all subgroup of the group $SF(\omega)$. 
Since the corresponding tower cardinal 
is involved in this question, the general 
problem of description the corresponding 
van Douwen's diagram for this lattice seems relevant. 
The paper is devoted to these questions.

\subsection{Almost containedness}

Van Douwen's diagrams are due to \cite{vd} and
\cite{vaughan}, where the case of $P(\omega)$ 
was considered. 
The term was used in \cite{c} where 
the case of partitions was studied. 
The general idea can be described as follows. \parskip0pt

Let $L$ be a lattice with 0 and 1, and let $I$ be an ideal of $L$.
We say that $a,b \in L \setminus I$ are 
{\it orthogonal} if $a \wedge b \in I$. 
The element $a$ is {\it almost contained} in $b$
(we denote it by $a \le_{a} b$) 
if $a \le b \vee c$ for some $c \in I$.
We write $a =_{a} b$ if $a \le_{a} b$ and $b \le_{a} a$. 
For any $a \in L$ we put $a_{I} = \{ b : b =_{a} a \}$. 
It is clear that the relation $\le_{a}$ becomes 
the usual almost containedness if we
consider the lattice $(P(\omega),\subseteq )$ 
with respect to the ideal
of finite subsets of $\omega$. \parskip0pt

In general, to characterize a lattice $L$ under these 
relations we need some further notions. 
We say that $a$ {\it splits} $b$
if there are $c,d \le b$ not in $I$ such that $c \le a$ and
$d,a$ are orthogonal.  
A family $\Gamma \subset L \setminus I$ is a
{\it splitting family} if for every $b \in L \setminus I$ 
there exists $a \in \Gamma$ that splits $b$. 
We say that $\Gamma$ is a {\it reaping family} if 
for each $a \in L \setminus I$ there is some 
$b \in \Gamma$ such that $b \le_{a} a$ or $a,b$
are orthogonal. 
We also define a family
$\Gamma \subset L \setminus I$ to be $\le$-{\em centered} 
if any finite intersection of its elements is not in $I$. \parskip0pt

We can now associate to $L$ the following cardinals. 
Define ${\mathfrak a}_{I}$ to be the least cardinality 
of an infinite  maximal family of pairwise orthogonal elements 
from $L \setminus 1_{I}$. 
Let ${\mathfrak p}_{I}$ be the least cardinality 
of a $\le$-centered family $\Gamma$ such that there 
is no $b \in L \setminus I$ such that $b$
is a lower bound of $\Gamma$ under $\le_{a}$ and the family
$\Gamma \cup \{ b \}$ is still $\le$-centered. 
Similarly, define ${\mathfrak t}_{I}$ (the tower cardinal) 
as the least cardinality of a $\le_{a}$-decreasing 
$\le$-centered chain without lower $\le_{a}$-bound
consistent (in the sense of $\le$-centeredness) with the family.
The cardinals ${\mathfrak s}_{I},{\mathfrak r}_{I}$ are the corresponding
(least) cardinals for splitting families and reaping
families respectively. 
It is worth noting that ${\mathfrak p}_{I}$ and
${\mathfrak t}_{I}$ can be undefined (for example if for any $a \in L$
the set $\{ b : b \le a \}$ is finite). 
Also, $(L,I)$ does not necessarily have a splitting family 
(for example if $L$ is an atomic 
boolean algebra and $I$ is trivial).
So ${\mathfrak s}_{I}$ can be undefined too. 
On the other hand, it is clear that 
${\mathfrak p}_{I} \le {\bf t}_{I}$ if they are defined. \parskip0pt

The last cardinal ${\mathfrak h}_{I}$ is defined as follows. A family $\Sigma$
of maximal families of pairwise orthogonal elements in
$L \setminus 1_{I}$ is {\it shattering} if for every
$a \in L \setminus I$ there are $\Gamma \in \Sigma$ and
distinct $b,c \in \Gamma$ which are not orthogonal to $a$.
Let ${\mathfrak h}_{I}$ be the least cardinality of a shattering family in $L$.
\parskip0pt

The following lemma seems to be folklore.

\begin{LEM} \label{1_1}
If ${\mathfrak s}_{I}$ is defined, then ${\mathfrak h}_{I} \le {\mathfrak s}_{I}$.
\end{LEM}

{\em Proof.}
Take a splitting family
$\Gamma = \{ c_{\nu} : \nu < {\mathfrak s} \}$. For each $\nu < {\mathfrak s}$
choose $\Psi_{\nu}$ to be a maximal family of pairwise orthogonal elements
such that $c_{\nu} \in \Psi_{\nu}$. Let us check that the set of
these families is shattering. Let $c \in L \setminus I$. Since
$\Gamma$ is a splitting family there is $\nu$ and $a,b \le c$ such
that $a \le c_{\nu}$ and $b$ is orthogonal to $c_{\nu}$. By our
construction there is $d \in \Psi_{\nu}$ not orthogonal to $b$.
So $\Psi_{\nu}$ shatters $c$ by $c_{\nu}$ and $d$. $\Box$

\bigskip

{\em Remark.} 
In the case of the lattice $(P(\omega), \cup , \cap )$
and the ideal $[\omega]^{<\omega}$ of all finite 
subsets of $\omega$ the introduced numbers are exactly 
the classical cardinals  $\mathfrak  a$, ${\mathfrak h}, {\mathfrak p}$, 
${\mathfrak r}, {\mathfrak s}, {\mathfrak t}$   
(all of them occur in \cite{vaughan}). 
Indeed, our definitions of 
${\mathfrak a}_{I}, {\mathfrak r}_{I}, {\mathfrak s}_{I}, {\mathfrak h}_{I}$
are formulated as the corresponding classical 
ones in \cite{vd} and in \cite{vaughan}.
The classical ${\mathfrak t}$ is the least cardinality 
of a $\le_{a}$-decreasing chain in $P(\omega)$ 
without $\le_{a}$-bound.
The classical ${\mathfrak p}$ is defined as follows.
We say that a family $\Gamma \subseteq [\omega]^{\omega}$
is $\le_{a}$-centered if every finite  $\Gamma'\subseteq \Gamma$ has an infinite pseudointersection: 
a set $X \in [\omega]^{\omega}$ almost contained in
each element of $\Gamma'$. 
Then the classical ${\mathfrak p}$ is the least
cardinality of a $\le_{a}$-centered family from 
$P(\omega)$ without lower $\le_{a}$-bound.   
So, there is no assumption on $\subseteq$-centeredness 
as in the definitions of ${\mathfrak p}_{I}$ and
${\mathfrak t}_{I}$. 
On the other hand we do not need such assumptions because
any $\le_{a}$-centered family from $P(\omega)$ is centered.
So ${\mathfrak p} = {\mathfrak p}_{[\omega]^{<\omega}}$ and
${\mathfrak t} = {\mathfrak t}_{[\omega]^{<\omega}}$. $\Box$

\bigskip

Note that the definitions of the above cardinals make sense if we consider
$L/=_{a}$ under the reverse order $\ge_{a}$ replacing the ideal $I$ by
$1_{I}$. In the case of $P(\omega)$ the converse cardinals are equal
to the corresponding cardinals for $\subset$ because $P(\omega)$ is a
Boolean algebra. The fact that this is not true in general is quite
important for the lattice of subgroups of $SF(\omega)$ and for partitions.
\parskip0pt

In the latter case we can consider partitions 
as subsets of $\omega^{2}$ under the inclusion 
(denoted by $\subset_{pairs}$).
The lattice that we get 
(with operations $\vee_{pairs}$ and $\wedge_{pairs}$) 
is converse to the lattice $\LATTICE$. 
Let $IF$ be the ideal of partitions (in $(\PART, <_{pairs})$) 
obtained from $id_{\omega}$ by adding a finite set of pairs. 
Then the class ${\mathfrak 1}_{IF}$ is exactly $\FPART$. 
Note that cardinal invariants of this lattice are studied in papers \cite{ma3} and \cite{brendle2}.
On the other hand, the relation of
almost containedness of partitions analyzed
in \cite{mat} and \cite{c} can be defined as follows:
$$
Y \le^{*} X \leftrightarrow (\exists Z \in IF )
(X \subset_{pairs} Y \vee_{pairs} Z).
$$
As a result we see that the cardinal 
${\mathfrak a}_{d},{\mathfrak p}_{d},{\mathfrak t}_{d}$, 
${\mathfrak h}_{d},{\mathfrak s}_{d},{\mathfrak r}_{d}$ 
examined in \cite{mat} and \cite{c} are 
the converse cardinals
for the pair $((\PART, <_{pairs}), IF)$.

\subsection{The lattice of subgroups of $SF(\omega)$}

Let $SF(\omega)$ be the group of all finitary 
permutations of $\omega$. 
This means that the elements of $SF(\omega)$ 
are exactly the permutations $g$ with finite support, 
where $supp(g) = \{ x: g(x) \not= x \}$. 
The algebraic structure of subgroups of $SF(\omega)$ 
is described in \cite{pn1}, \cite{pn2}. 
The aim of our paper is to study the van Douwen's 
invariants of the lattice of
subgroups of $SF(\omega)$. \parskip0pt

Throughout the paper $LF$ is the lattice of all subgroups of
$SF(\omega)$ and $IF$ is the ideal of all finite subgroups.
We say that $G_{1}$ and $G_{2}$ from $LF \setminus IF$ are {\em orthogonal}
if their intersection is in $IF$. The group $G_{1}$ is {\em almost contained}
in $G_{2}$ ($G_{1} \le_{a} G_{2}$) if $G_{1}$ is a subgroup of a group
finitely generated over $G_{2}$ by elements of $SF(\omega)$.
Let $SF(\omega)_{IF} = \{ G \le SF(\omega): SF(\omega)$ is finitely
generated over $G \}$. As in Section 1.2 we define the cardinal numbers
${\mathfrak a}_{SF}$, ${\mathfrak p}_{SF}$, ${\mathfrak t}_{SF}$, ${\mathfrak r}_{SF}$,
${\mathfrak h}_{SF}$ and ${\mathfrak s}_{SF}$. For example,
${\mathfrak a}_{SF}$ is the least cardinality of a maximal family of
pairwise orthogonal elements from $LF \setminus SF(\omega)_{IF}$ and
${\mathfrak p}_{SF}$ is the least cardinality of a $\le$-centered
family of elements in $LF \setminus IF$ with no lower $\le_{a}$-bound
$\le$-consistent (in the sense of $\le$-centeredness) with the family.

We put a topology on $LF$ in the following way. Let $H \le SF(\omega)$
be a finite group and $A \subset \omega$ be a finite set containing the union of
the supports of the elements of $H$. Let $[H,A]$ be the set of all
subgroups of $SF(\omega)$ such that the groups  they induce on
$A$ are equal to $H$ (we think of $H$ as a permutation group on $A$).
The topology that we consider is defined by the base consisting of all
sets $[H,A]$. This topology is metrizable: fix an enumeration
$A_{0},A_{1},...$ of all finite subsets of $\omega$ and define

\begin{quote}
$d(G_{1},G_{2}) = \sum \{ 2^{-n}:$ the groups induced by
$G_{1}$ and $G_{2}$ on $A_{n}$ are not the same~$\}$.
\end{quote}

Note that the space $LF$ is complete. A function
$\delta: LF \rightarrow n, n \in \omega,$ is then called a
Borel (respectively ${\bf \Sigma}^{1}_{1} \cup {\bf \Pi}^{1}_{1}$)
coloring if $\delta^{-1}(i)$ is Borel (respectively analytic or
coanalytic) for every $i < n$ 
(where $n \in \omega$ is viewed as $\{ 0, . . . ,n - 1 \}$.)  
\parskip0pt

Consider the set $LF_{1}$
of all groups of the form $SF(\omega) \cap G$ where $G$ is 1-closed.
We identify elements of $LF_{1}$ with elements of
$2^{\omega \times \omega}$ (the corresponding partitions). 
Then
it is easily seen that the topology on $LF_{1}$ induced by the topology
above becomes the restriction of the product topology on
$2^{\omega \times \omega}$ where $2$ is considered discrete.
A theorem of T. Carlson and S. Simpson from \cite{b} can be restated as
follows: for every
$({\bf \Sigma}^{1}_{1} \cup {\bf \Pi}^{1}_{1})$-coloring
$\delta: LF_{1} \rightarrow 2$ there exists
$G \in LF_{1} \setminus SF(\omega)_{IF}$ such that $\delta$
is constant on the elements of $LF_{1} \setminus SF(\omega)_{IF}$
containing $G$. The corresponding theorem for $P(\omega)$ proved by
F. Galvin and K. Prikry in \cite{GalPrik} is stated as follows:
for every  $({\bf \Sigma}^{1}_{1} \cup {\bf \Pi}^{1}_{1})$-coloring
(originally: Borel coloring; see Remark 2.6 in \cite{b})
$\delta: P(\omega) \rightarrow 2$ there exists an infinite
$A \in P(\omega)$ such that $\delta$ is constant on the set of all
infinite subsets of $A$. It is shown in \cite{b} that this theorem is
a consequence of the Carlson-Simpson theorem.

\section{The diagram in $(LF,IF)$}

In this section we describe the diagram for the lattice $(LF, IF)$. 
First we compare the coefficients of this diagram with their classical analogues.

\subsection{Comparing $(LF,IF)$ with $(\mathcal{P}(\omega ), Fin)$} 

The following  result specifies basic relations between the coefficients for $(LF, IF)$ and  the coefficients for 
$(P(\omega)/fin, \subseteq_*)$.\bigskip

\begin{PRO} \label{srph} 
The following inequalities are true in {\bf ZFC}.
$${\mathfrak t}_{SF}\le {\mathfrak t},\  {\mathfrak p}_{SF}\le  {\mathfrak p},\  {\mathfrak r}_{SF} \le {\mathfrak r},\  
{\mathfrak s}\le {\mathfrak s}_{SF}\ \mbox{ and }{\mathfrak h}\le {\mathfrak h}_{SF}.$$

\end{PRO} 

{\em Proof.} 
Fix an arbitrary group  $G^*\subseteq SF(\omega)$  generated by the set
$\{\sigma_i: i \in \omega\}$ of pairwise disjoint cycles  with consequtive prime orders, 
 i.e.  $|supp (\sigma_i)| =p_i$, for $i\in\omega$.
Such a  $G^*$ is isomorphic to the direct sum of 
all $\mathbb{Z}/p_{i}\mathbb{Z}$, $i\in\omega$.
Denote by  $LF_{G^*}\le LF$ the sublattice of all elements below $G^*$.
To each $A\subset \omega$ we associate 
the subgroup  $G_A<G^*$ generated by 
all $\sigma_i$ with $i\in A$. 
It is clear that the map $A\to G_A$ induces an isomorphism from $(\mathcal{P}(\omega) /Fin ,\  \subseteq_{a} )$ 
to  $(LF_{G^*}/=_a, \le_a) $.

Now the first two inequalities are obvious.

To see that ${\mathfrak r}_{SF}\le {\mathfrak r}$, observe that for any family ${\mathcal R}\subseteq LF_{G^*}$ reaping for $LF_{G^*}$, 
the family ${\mathcal R}\cup \{{G^*}\}$ is reaping for $LF$. 

To prove ${\mathfrak s}\le {\mathfrak s}_{SF}$, note that for any  family  
${\mathcal S}\subseteq LF$  splitting for $LF$,  the family $\{K\cap  {G^*}: K\in {\mathcal S}\}\setminus IF$ is splitting for   $LF_{G^*}$.

For ${\mathfrak h}\le {\mathfrak h}_{SF}$ we argue as follows.
 If ${\mathcal A}$ is a maximal family  of  pairwise orthogonal elements from  $(LF\setminus SF(\omega)_{ IF})$,
then the  family ${\mathcal A}_{G^*}=\{K\cap  {G^*}: K\in {\mathcal A}\}\setminus IF$  
 is  maximal for $LF_{G^*}$,  although it is not necessarily infinite.
Nevertheless, for any family $\mathcal{H}$ shattering  for $LF$, 
the family ${\mathcal H}_{G^*}=\{{\mathcal A}_{G^*}: {\mathcal A}\in {\mathcal H}, {|\mathcal A}_{G^*}|>1\}$ 
is  nonempty.  

For  every finite ${\mathcal B}\in{\mathcal H}_{G^*}$, where ${\mathcal B}=\{H_0, H_1, \ldots , H_{n-1}, H_n\}$,
choose an arbitrary infinite maximal  pairwise orthogonal family 
${\mathcal B'}\subseteq LF_{G^*}$ such that $\{H_0, H_1, \ldots , H_{n-1}\}\subseteq {\mathcal B'}$.
It is easy to see that the family $$\{{\mathcal B}\in {\mathcal H}_{G^*}: {\mathcal B}\mbox{ is infinite}\}\cup
\{{\mathcal B'}: {\mathcal B}\in {\mathcal H}_{G^*}\mbox{ and } {\mathcal B}\mbox{ it  is finite}\}$$ is shattering for $LF_{G^*}$.
$\Box$ 
\bigskip 

The next result reduces the first two inequalities to the following equality.

\begin{PRO} \label{2_2} 
$$ 
{\mathfrak p}_{SF} = {\mathfrak p} = {\mathfrak t}_{SF}  \mbox{ . } 
$$
\end{PRO} 

{\em Proof.} Fix an arbitrary one-to-one enumeration  
$SF(\omega)=\{\rho_i: i\in\omega\}$.

We have to prove  ${\mathfrak p}\le {\mathfrak p}_{SF}$.
Take a $\le$-centered family 
$\Gamma \subset LF \setminus IF$ of cardinality 
${\mathfrak p}_{SF}$ without a $\le_a$- bound $H\in LF\setminus IF$ such that $\Gamma\cup H$  is still $\le$-centered.
Then to every $G\in \Gamma$ assign the set $A_G=\{i: \rho_i\in G\}$.
 It is obvious that  $\{A_G:G\in \Gamma\}$ is a $\subseteq$--centered family of subsets.
We claim that it does not have an  infinte $\subseteq^*$-bound $A$ such, that the family
$\{A_G:G\in \Gamma\}\cup \{A\}$ is $\subseteq$-centered.
Suppose the contrary and let $A$ be such a bound.
Then   the group $G_A$ generated by the set $\{\rho_i: i\in A\}$ is almost contained in every $G\in \Gamma$ and the family 
$\Gamma\cup \{G_A\}$ is $\le$-centered.
This contradicts the assumption.
Therefore, ${\mathfrak p} \le {\mathfrak  p}_{SF}$.  

Now it suffices to  use the recent theorem of Melliaris and Shelah that ${\mathfrak p}={\mathfrak t}$ (see  \cite{ms})  and
Proposition \ref{srph} to complete the proof.
$\Box$

\subsection{The diagram}
We need a couple of lemmas.

\begin{LEM}  \label{1-lemma} 
Let $G \in LF \setminus IF$ and $m \in \omega$. Then:
\begin{enumerate}
\item [{(i)}] there exists a non-trivial $\rho \in G$ such that
$supp(\rho) \cap m = \emptyset$ ,
\item [(ii)] moreover, for any $H \subset Sym(m)$ and any sequence
$G_{0},G_{1}, . . . ,G_{n} \in LF \setminus IF$ of groups orthogonal
to $G$ the above $\rho$ can be chosen so that additionally
$\langle H,\rho \rangle \cap G_{i} = \langle H \rangle \cap G_{i}$, for
$i \le n$.
\end{enumerate}
\end{LEM}

{\em Proof.  (i).}  Suppose that the lemma is not true. Choose
a minimal $A = \{ a_{0}, . . . , a_{k} \} \subset m$ such that there
are infinitely many $g \in G$ satisfying 
$m \cap supp(g) \subset A$. Then $A \not= \emptyset$. We fix some
non-trivial $g_{0}$ with that property and consider all tuples
$g(\bar{a}) = (g(a_{0}), . . . , g(a_{k}))$ for the above
$g$'s. If each of these tuples has non-empty intersection with
$supp(g_{0})$, then there is $i \le k$ such that $g(a_{i})$ is
the same for infinitely many $g$'s. Clearly, for such $g$ and
$g'$ the set $m \cap supp(g^{- 1} \cdot g')$ is
a subset of $A \setminus \{ a_{i} \}$. This contradicts the
minimality of $A$. \parskip0pt

Choose $g$ as above with
$g(\bar{a}) \cap supp(g_{0}) = \emptyset$ additionally.
It is easily seen that
$g^{- 1} \cdot g_{0} \cdot g$ fixes
$m$ pointwise. This contradicts our assumption. \parskip0pt

{\em (ii). } Suppose the contrary. By ({\em i}) we can find
$i \le n$ such that for infinitely many $\rho \in G$ with
$supp(\rho) \cap m = \emptyset$ there is $g \in \langle H \rangle$ satisfying
$g \cdot \rho \in G_{i}$. Since $\langle H \rangle$ is finite, there is
$g_{0} \in \langle H \rangle$ such that $g_{0} \cdot \rho \in G_{i}$
for infinitely many $\rho \in G$. Hence, for infinitely many
$\rho,\rho' \in G$, $\rho^{-1} \cdot \rho' \in G_{i}$,
which contradicts orthogonality. $\Box$

\bigskip

As a consequence of the above lemma we get the following easy statement 
which is a corollary of Proposition \ref{2_2} either.

\begin{LEM} \label{chain} 
For any countable sequence $G_{0} > G_{1} >...$ of elements of
$LF \setminus IF$ there is a group $G \in LF \setminus IF$ such that
$G \le_{a} G_{i}$, for every $i \in \omega$,  and the family
$\{ G_{i}: i \in \omega \} \cup \{ G \}$ is $\le$-centered.
\end{LEM}

{\em Proof. }
Assume  we have a decreasing sequence
$G_{0} > G_{1} > . . . $ in $LF \setminus IF$.
For every $i \in \omega$ choose non-trivial $g_{i} \in G_{i}$ such that
$supp(g_{i})$ is disjoint from the supports of the previous elements.
We can do this by Lemma \ref{1-lemma}(i). 
Let $G$ be the group generated by all
these $g_{i}$. 
Then $G \in LF \setminus IF$, the family
$\{ G \} \cup \{ G_{i} : i \in \omega \}$ is centered and
$G \le_{a} G_{i}$, for every $i \in \omega$. 
$\Box$ 

\bigskip

\begin{LEM} \label{3_1} 
Let $G_{0}, . . . ,G_{n - 1}$ be a sequence of infinite groups from
$LF$ not $a$-equivalent to $SF(\omega)$. Then for any $k,m \in \omega$,
$k > 0$, and $H \subset Sym(m)$ there is a non-trivial finitary permutation
$\rho$ consisting of $(k + 1)$-cycles such that
$supp(\rho) \subset \omega \setminus m$ and for every $i < n$,
$$
\langle H,\rho \rangle \cap G_{i} = \langle H \rangle \cap G_{i}.
$$
\end{LEM}

{\em Proof.}  For each $i < n$ set
$$
S_{i} = \{ g \in G_{i}:
\exists g_{0},g_{1} (g_{0} \in \langle H \rangle \wedge
(m \cap supp(g_{1}) = \emptyset) \wedge (g = g_{0} \cdot g_{1})) \}
$$
It is easily seen that each $S_{i}$ is a group.
Choose a family $\{ D_{j0}: 0 \le j < n \}$ of pairwise disjoint
finite sets such that for every $j$, $D_{j0} \subset \omega \setminus m$
and $S_{j}$ does not induce $Sym(D_{j0})$.
Let $D_{j1},. . . ,D_{jk}$ be sets from $\omega \setminus m$ of the
same size as $D_{j0}$. We may assume that every pair from
$\{ D_{ji} : 0 \le i \le k; j < n \}$ has empty intersection.
For every $0 < i \le k$ and $j < n$ we choose a bijection
$f_{ji}$ from $D_{j0}$ onto $D_{ji}$ such that it is not induced by
any element of $S_{j}$. The existence of such $f_{ji}$ is a consequence of
the fact that for any bijections $f,g : D_{j0} \rightarrow D_{ji}$
induced by $S_{j}$ , the bijection $g^{- 1} \cdot f$ defines
a permutation on $D_{j0}$ induced by $S_{j}$.
\parskip0pt

We now define a permutation $\rho$ with the support
$\bigcup \{ D_{ji} : 0 \le i \le k, 0 \le j < n \}$ as follows.
If $x \in D_{ji}, 0 < i < k,$ then $\rho(x) = f_{j(i + 1)}(f^{-1}_{ji}(x))$.
If $x \in D_{j0}$, then $\rho(x) = f_{j1}(x)$. For $x \in D_{jk}$ we put
$\rho(x) = f^{-1}_{jk}(x)$. Let us check
that $\rho$ satisfies the conclusion of the lemma. It is clear that
$\rho$ consists of cycles of length $k + 1$. Suppose, that
for some $g_{0} \in \langle H \rangle$ the element
$g = g_{0} \cdot \rho^{l}, 0 < l \le k,$
is contained in some $G_{j}$. Thus $g \in S_{j}$ and by our construction
$g$ maps $D_{j0}$ onto $D_{jl}$ by $f_{jl}$.
Since $S_{j}$ does not induce $f_{jl}$,
we have a contradiction. $\Box$

\bigskip

The van Douwen cardinals for
$(LF,IF)$ are described in the following theorem.

\begin{THM}\label{diagram}
\begin{itemize}
\item[(i)] The following inequalities are true in $(LF,IF)$:
$$
\omega_{1} \le {\mathfrak p}_{SF} = {\mathfrak t}_{SF} \le {\mathfrak h}_{SF}
\le {\mathfrak s}_{SF} \le 2^{\omega},
$$
$$
\omega_{1} \le {\mathfrak a}_{SF}, {\mathfrak r}_{SF} \le 2^{\omega};
$$

\item[(ii)] All the coefficients are equal to continuum under
Martin's Axiom;

\item[(iii)] Each of the following equalities is consistent with
$\{ {\bf ZFC} + \omega_{1} < 2^{\omega} \}$:
$$
{\mathfrak a}_{SF} = \omega_{1}, {\mathfrak s}_{SF} = \omega_{1},
{\mathfrak r}_{SF} = \omega_{1}.
$$

\end{itemize}
\end{THM}

{\em Proof. (i).} The inequality ${\mathfrak h}_{SF} \le {\mathfrak s}_{SF}$
is shown in Lemma \ref{1_1}. 
The inequality $\omega_{1} \le {\mathfrak p}_{SF}$
follows from Proposition \ref{2_2} and
${\mathfrak t}_{SF} \le {\mathfrak h}_{SF}$    is a consequence of Propositions \ref{srph} and  \ref{2_2} 
together with the classical inequality ${\mathfrak t}\le {\mathfrak h}$. 
 \parskip0pt

To prove $\omega_{1} \le {\mathfrak r}_{SF}$ it suffices to show that
if a family $\Psi \subseteq LF \setminus IF$ is countable then
there exists $G \in LF \setminus IF$ such that for every
$G' \in \Psi$ the groups $G, G'$ are not orthogonal and
$G' \not\le_{a} G$. 
Let $\{ G_{0}, G_{1}, . . . \}$ be an enumeration
of $\Psi$. 
Assume that each member of $\Psi$ occurs infinitely often.
We construct two sequences $g_{0}, g_{1}, ...$
and $h_{0}, h_{1}, ... $ of finitary permutations with pairwise
disjoint supports such that for all $i,j \in \omega$ we have
$supp(g_{i}) \cap supp(h_{j}) = \emptyset$ and
$g_{i}, h_{i} \in G_{i}$. 
It is easily seen that Lemma \ref{1-lemma}(i) implies
the existence of such sequences. 
Let $\hat{G}_{1} = \langle \{ g_{i}: i \in \omega \} \rangle$
and $\hat{G}_{2} = \langle \{ h_{i}: i \in \omega \} \rangle$. 
Clearly, $\hat{G}_{1}$ and $\hat{G}_{2}$ are orthogonal 
but they are not orthogonal to any $G_{i}$
(since each member of $\Psi$ is enumerated infinitely often). 
Now it is easy to see that $G = \hat{G}_{1}$ 
satisfies the conditions that we need. \parskip0pt

To prove the inequality $\omega_{1} \le {\mathfrak a}_{SF}$ take 
a countable $\Psi \subset LF \setminus SF(\omega)_{IF}$. 
We construct a group $G$ by induction. 
Fix an enumeration of $\Psi : G_{0},G_{1},...$ .
Let $H$ be the set of the elements which 
have been constructed at the first $n - 1$ steps.
At the $n$-th step we choose a permutation $\rho$ as in 
Lemma \ref{3_1} with respect to $G_{0},...,G_{n}$ and $m$ large enough. 
It is easily seen that that the group generated by 
this sequence is orthogonal to any group from $\Psi$. \parskip0pt

\medskip

{\em (ii).} Assume {\bf MA}. By Proposition \ref{srph},  \ref{2_2}  and the classical result 
${\bf MA}\models {\mathfrak p}=2^{\omega}$ we have
$${\mathfrak p}_{SF} ={\mathfrak  t}_{SF}={\mathfrak s}_{SF}={\mathfrak h}_{SF}= 2^{\omega}.$$

To prove ${\mathfrak r}_{SF} = 2^{\omega}$ we introduce a $ccc$ forcing
notion ${\bf P}_{r}$ as follows. 
Consider the family of all pairs $(H,H')$
where $H,H' \subset SF(\omega)$ are finite and the supports of any two
elements of $H \cup H'$ have empty intersection. 
The order is defined as follows $(H, H') \le (F,F')$ iff 
$F \subseteq H$ and $F' \subseteq H'$.
Let $\Psi \subset LF \setminus IF$ have cardinality $< 2^{\omega}$.
For any $k \in \omega$ and $G \in \Psi$ the family
$$
\{ (H,H') \in {\bf P}_{r}: k < |H' \cap G|, k < |H \cap G| \} 
$$ 
is dense in ${\bf P}_{r}$ by Lemma \ref{1-lemma}(i) 
(see also the previous part of the proof). 
For a generic $\Phi$ define
$G_{0} = \langle \bigcup \{ H: (H,H') \in \Phi \} \rangle$. 
It is easy to see that for any $G \in \Psi$, the groups $G$ 
and $G_{0}$ are not orthogonal and $G$ is not contained 
in $G_{0}$ under $\le_{a}$. 
Thus $\Psi$ is not reaping. \parskip0pt

To show ${\mathfrak a}_{SF} = 2^{\omega}$, given an infinite
family $\Gamma \subset LF$ of infinite groups define 
a forcing notion ${\bf P}_{a}$ as follows. 
Let ${\bf P}_{a}$ be the set of all pairs
$(H,F)$ where $F$ is a finite subset of $\Gamma$ 
and $H$ is a finite set of permutations such 
that their supports are pairwise disjoint. 
We define $(H,F) \le (H',F')$ iff 
$H' \subset H, F' \subset F$ and each 
$h \in \langle H \rangle \setminus \langle H' \rangle$ 
is not contained in any $G \in F'$. 
It is easily verified that
${\bf P}_{a}$ is a {\it ccc} forcing notion. \parskip0pt

Consider ${\bf P}_{a}$ with respect to 
$\Psi \subset LF \setminus SF(\omega)_{IF}$ 
of cardinality $< 2^{\omega}$.
Clearly, the following sets are dense in 
${\bf P}_{a}$ (apply Lemma \ref{3_1} in the second case):
$$
\Sigma_{G} = \{ (H,F) : G \in F \}, G \in \Psi \mbox{ , and } 
$$ 
$$
\Sigma_{l} = \{ (H,F) : 
\mbox{ the number of the elements of } H 
\mbox{ is greater than } l \}, l \in \omega .
$$  
By ${\bf MA}$ we have a filter $\Phi \subset {\bf P}_{a}$
meeting all these $\Sigma$'s. It is easy to see that the group
$G_{0} = \langle \bigcup \{ H : (H,F) \in \Phi \} \rangle$ 
is orthogonal to any group from $\Psi$.

\medskip

{\em (iii).} It follows   from {\em (i)} and Proposition \ref{srph} that
$$Con\left({\bf ZFC} + ({\mathfrak p}_{SF} ={\mathfrak t}_{SF}={\mathfrak r}_{SF}= \omega_{1} < 2^{\omega})\right).$$
\medskip

To prove $Con({\bf ZFC} + {\mathfrak a}_{SF} = \omega_{1} < 2^{\omega})$.
we start with an arbitrary countable family
$\Psi_{0} \subset LF \setminus SF(\omega)_{IF}$
of parwise orthogonal groups. 
Take a sequence
$$
\Psi_{0} \subset \Psi_{1} \subset ... \subset \Psi_{\gamma} \subset ...,
\gamma < \omega_{1},
$$
by a finite support iteration
$$
({\bf P}_{\gamma}, Q_{\gamma}: \gamma < \omega_{1})
$$
of the forcing ${\bf P}_{a}$ (from the previous part of the proof) 
applied to potential $\Psi_{\gamma}$'s.  
The canonical name for ${\bf P}_{\gamma}$
of $\Psi_{\gamma + 1}$ is obtained from the canonical name
of $\Psi_{\gamma}$ by adding the canonical name of the
group $G_{\gamma}$ defined by $Q_{\gamma}$ as $G_{0}$
by ${\bf P}_{a}$ above. Let $\Phi$ be generic
for ${\bf P}_{\omega_{1}}$ and $\Phi_{\gamma}$ be the corresponding
restriction to ${\bf P}_{\gamma}$. 
It is easily seen that
${\bf P}_{\omega_{1}}$ fulfils the {\it ccc}. 
Since $\omega_{1}$ is regular and each group $G$ in
$LF[\Phi]$ is defined by a countable set
of finitary permutations,  it is contained
in some $LF[\Phi_{\gamma}]$. 
Suppose that some $G$ is orthogonal to each group from 
$\Psi_{\gamma}$. 
Thus by Lemma \ref{1-lemma}(ii) every set
$$
D_{n} = \{ (H,F) \in Q_{\gamma}[\Phi_{\gamma}] : n < |H \cap G| \}
$$
is dense in $Q_{\gamma}[\Phi_{\gamma}]$. 
So the group
$$
G_{\gamma} =  \big\langle \bigcup \{ H: (H,F) \in
\Phi_{\gamma + 1}/ \Phi_{\gamma} \} \big\rangle
$$
is not orthogonal to $G$. 
This shows that the set 
$\bigcup \{ \Psi_{\gamma}: \gamma < \omega_{1} \}$ 
is a maximal family of pairwise orthogonal groups in $LF[\Phi]$. \parskip0pt

The case $Con({\bf ZFC} + {\mathfrak  s}_{SF} = \omega_{1} < 2^{\omega})$
can be handled in a similar way - constructing $G_{\gamma}$
we apply the forcing ${\bf P}_{r}$. 
$\Box$

\bigskip

We conjecture that ${\mathfrak h}_{SF} = {\mathfrak h}$, 
${\mathfrak s}_{SF} = {\mathfrak s}$ and ${\mathfrak r}_{SF} = {\mathfrak r}$. 
Note that the corresponding equalities 
hold for the lattice of partitions under 
$\le_{pairs}$ \cite{brendle2}. 
At the moment we cannot adapt the arguments 
of \cite{brendle2} to our case. 
The case of ${\mathfrak a}$ is also open. 

\bigskip

\section{Two variants of Matet's theorem}

\subsection{ The first version of Matet's theorem} 

The following theorem is  formulated for the context described in the previuos section.

The proof of the theorem of Matet stated in Section 1.1  
(this is Proposition 8.1 from \cite{mat}) uses 
the Carlson-Simpson's theorem and Proposition 4.2 from \cite{mat} 
asserting that ${\mathfrak t}_{d}$ is uncountable. 
We will use the same strategy.

\begin{THM}
Assuming ${\bf MA}$ there is a filter 
$F \subset LF \setminus IF$ such that  for every
$({\bf \Sigma}^{1}_{1} \cup {\bf \Pi}^{1}_{1})$-coloring
$\delta : LF \rightarrow 2$ there is $G \in F$ 
such that $\delta$ is constant on the set of 
all infinite subgroups of $G$.
\end{THM}

{\em Proof.}
 Let $G^*$ be a  group generated by an infinite family 
$\{\sigma_i: i\in\omega\}$ of finite permutations with pairwise disjoint supports and distinct prime orders.
For example we can take the group described in the proof of Proposition  \ref{srph}. 
Then every $G \le {G^*}$  is generated by a subset of the set $\{ \sigma_i: i \in \omega \}$ and
 the lattice of all subgroups of ${G^*}$ 
is isomorphic to $(P(\omega), \subseteq)$.
We identify $G \le {G^*}$ with the corresponding subset of $\omega$.
Notice that then the topology defined in Section 1.3,  
on $\{ G : G \le {G^*} \}$ 
becomes the product topology on $2^{\omega}$. 
Also, $\{ G : G \le {G^*} \}$
is a closed subset of $LF$. \parskip0pt

We now use the strategy of Proposition 8.1 from \cite{mat}.
Let $\langle\delta_{\alpha}: \alpha < 2^{\omega}\rangle$ be an enumeration of all
$({\bf \Sigma}^{1}_{1} \cup {\bf \Pi}^{1}_{1})$-colorings
$\delta : LF \rightarrow 2$. We construct a
descending tower of subgroups of ${G^*}$.  Supposing that
$G_{\beta}, \beta < \alpha,$ have already been selected,
use Theorem \ref{diagram} (ii) to find $G_{\alpha} \le {G^*}$ such that
the family $\{ G_{\gamma} : \gamma \le \alpha \}$ is $\le$-centered
and $G_{\alpha} \le_{a} G_{\gamma}$ for all $\gamma < \alpha$.
By the Galvin-Prikry theorem (\cite{GalPrik}) there is an infinite subset
of the set of generators of $G_{\alpha}$ such that all its infinite
subsets have the same color with respect to the coloring induced by
$\delta_{\alpha}$. This shows that $G_{\alpha}$ can be chosen such that
all its infinite subgroups have the same color with respect to
$\delta_{\alpha}$. \parskip0pt

Let $F$ be the filter generated by the tower obtained. It follows from
the construction that $F$ satisfies the conditions of the theorem.
$\Box$

\subsection{Another version of Matet's theorem}

As we noted in Introduction the lattice of partitions under
the reverse order is a sublattice of $LF$. 
This suggests that in the lattice $LF$ the most natural variant 
of the theorem of P. Matet cited there 
(Proposition 8.1 from \cite{mat}) is the folowing one.

\begin{THM} \label{matth}
Assuming the continnuum hypothesis there is an ideal
$I \subset LF \setminus SF(\omega)_{IF}$
such that  for every
$({\bf \Sigma}^{1}_{1} \cup {\bf \Pi}^{1}_{1})$-coloring
$\delta : LF \rightarrow 2$ there is $G \in I$ such that
$\delta$ is constant on the set of all supergroups of $G$
which do not belong to $SF(\omega)_{IF}$.
\end{THM}

In the proof of the statement  we shall apply the result below.

\begin{LEM} \label{super} 
Let $P_{1},...,P_{i},...$ be a sequence of pairwise 
disjoint infinite subsets of $\omega$ defining a partition 
$E_{0}$ of $\omega$.
Let 
$G_{0} = Aut(\omega, P_{1},...,P_{i},...) \cap SF(\omega)$ 
be the subgroup of $SF(\omega)$ corresponding 
to a 1-closed subgroup of $Sym(\omega)$ defined 
by $P_i$, $i\in \omega\setminus \{ 0\}$. 
Then any proper supergroup of $G_{0}$
has this form for a partition coarser than $E_{0}$. 
\end{LEM} 

{\em Proof.} 
Let $g$ be a finitary permutation such that 
$g(a) = b \in P_{j}$ for $a \in P_{i}, i \not= j$. 
Let $a' \in P_{i} \setminus supp(g)$ and
$b' \in P_{j} \setminus supp(g)$. 
Below we denote the transposition
of $x$ and $y$ by $(x,y)$. 
It is clear that the element
$(a,a')\cdot g^{-1}\cdot (b,b')\cdot g\cdot (a,a')$ 
(which belongs to $\langle G_{0},g \rangle$) 
is the transposition $(a',b')$. 
This yields that the group inducing $SF(P_{i} \cup P_{j})$ 
and acting trivially on $\omega \setminus (P_{i} \cup P_{j})$, 
is a subgroup of $\langle G_{0},g \rangle$. 
The rest is clear. 
$\Box$ \bigskip

{\em Proof of the theorem.}
Let $P_{1},...,P_{i},...$ be a sequence of paiwise disjoint infinite
subsets of $\omega$ defining a partition $E_{0}$ of $\omega$.
Then $G_{0} = Aut(\omega, P_{1},...,P_{i},...) \cap SF(\omega)$ is
the subgroup of $SF(\omega)$ obtained from the corresponding 1-closed
subgroup of $Sym(\omega)$. 
By Lemma \ref{super} any proper supergroup of $G_{0}$
has this form for a partition coarser than $E_{0}$. \parskip0pt

We may now consider the set 
$L_{0} = \{ G: G_{0} \le G \le SF(\omega) \}$
as a sublattice of partitions coarser than $E_{0}$. 
Notice that then the topology defined in Section 1.3, 
on $L_{0}$ becomes the product topology on
$2^{\omega \times \omega}$. 
This follows from the fact that any finite
permutation group (on a finite subset of $\omega$) 
induced by a group $G$ from $L_{0}$ is a finite 1-closed 
permutation group and can be identified with a partition 
induced by the partition corresponding to $G$. 
Moreover, it is easy to see that 
$L_{0}$ is closed in $LF$. \parskip0pt

We now use the Matet's theorem. 
Take an ideal $I_{0}$ of
$L_{0}$ provided by this theorem. 
Then $I_{0}$ generates an ideal of $LF$. 
This ideal works as $I$ in the statement.
$\Box$

\section{Remarks}

\subsection{The dual diagram for $LF$} 

Theorem \ref{matth} suggests investigation of 
the reverse ordering of $LF$. 
Using  Lemma \ref{super} we get a result analogous to Proposition \ref{srph}.

\begin{PRO} 
Let 
${\mathfrak h}_d, {\mathfrak p}_d, {\mathfrak r}_d, {\mathfrak s}_d, {\mathfrak t}_d$ 
be dual cardinal invariants of the lattice of partitions 
defined by the scheme of Section 1.2 (defined as in \cite{c}). 
Let us consider $LF /=_a$ with respect to the converse 
ordering $\ge_a$ and let 
${\mathfrak h}^d_{SF}$, ${\mathfrak p}^d_{SF}$,${\mathfrak r}^d_{SF}$,${\mathfrak s}^d_{SF}$, ${\mathfrak t}^d_{SF}$ 
be the corresponding sequence of cardinal invariants 
defined with respect to $SF(\omega )_{IF}$ as 
an ideal of this converse lattice. 

Then  
${\mathfrak h}_d \le {\mathfrak h}^d_{SF}$,  ${\mathfrak r}^d_{SF}\le {\mathfrak r}_d$, 
${\mathfrak s}_d \le {\mathfrak s}^d_{SF}$ 
and 
 ${\mathfrak t}^d_{SF} = {\mathfrak p}^d_{SF} = \omega_1$. 
\end{PRO} 

{\em Proof.} 
We take  any one-closed group defined by a partition into infinitely many classes and use 
Lemma \ref{super} to argue as in the proof of  Proposition  \ref{srph} to obtain
$$
{\mathfrak h}_d \le {\mathfrak h}^d_{SF} \mbox{ ,  } {\mathfrak r}^d_{SF}\le {\mathfrak r}_d \mbox{ , } 
{\mathfrak s}_d \le {\mathfrak s}^d_{SF} \mbox{ , } 
{\mathfrak t}^d_{SF}\le {\mathfrak t}_d \mbox{ , } {\mathfrak p}^d_{SF}\le {\mathfrak p}_d . 
$$ 
Since ${\mathfrak p}_d = {\mathfrak t}^d_{SF} = \omega_1$ (see \cite{mat}), 
we have the statement of the lemma. 
$\Box$ 

\bigskip

Using this proposition and the material of 
papers \cite{mat}, \cite{c}, \cite{brendle1}, \cite{minami} 
we obtain the following relations: 
$$ 
{\mathfrak r}^d_{SF}\le {\mathfrak r}_d \le min ({\mathfrak r}, {\bf {\mathfrak d}}, non(\mathcal{M}), non(\mathcal{N})) 
\mbox{ and } 
max ( cov(\mathcal{N}), cov (\mathcal{M}), {\mathfrak s}, {\mathfrak b}) \le {\mathfrak s}_d \le {\mathfrak s}^d_{SF} . 
$$
Moreover the following relations are consistent with {\bf ZFC}: 
$$
{\mathfrak r}_{d}\le add(\mathcal{M})\mbox{ , } {\mathfrak r}_d > {\mathfrak b} \mbox{ , } 
{\mathfrak s}_d > cof (\mathcal{M}) \mbox{ , } {\mathfrak s}_d \le {\mathfrak r} \mbox{ , } {\mathfrak s}_d < {\bf {\mathfrak d} }. 
$$
We mention the following 
questions: 

1. Is ${\mathfrak a}^d_{SF} = 2^{\omega}$? 

2. What is relation between ${\mathfrak h}^d_{SF}$ and ${\mathfrak h}$? 

3. Are the following relations consistent with {\bf ZFC}: 
$$
{\mathfrak r}^{d}_{SF} > {\mathfrak b} \mbox{ , } 
{\mathfrak s}^d_{SF} \le {\mathfrak r} \mbox{ , } {\mathfrak s}^d_{SF} < {\mathfrak {\mathfrak d} }?  
$$

\subsection{Other remarks}

Our results suggest the investigation
of van Douwen's cardinals for the lattice of all closed subgroups of
$Sym(\omega)$. The definition of the $a$-order in this case must be
as follows: $G \le_{a} G'$ iff there exists a finite
set $X$ of finitary permutations such that $G$ is a subgroup
of the closed group generated by $G'$ and $X$. It is worth noting that
some results of \cite{mac} can be interpreted in this vein for some
converse coefficients (for example, see Observation 3.3 in \cite{mac}).
\parskip0pt

However, one can notice that the lattice of all closed subgroups
admits several constructions which make some of the van Douwen's cardinals
trivial. 
For example, the group $\mathbb{Z}$ with the natural action on itself
can be considered as a closed subgroup of $Sym(\omega)$. 
It is
clear that for every $n \in \omega$ no closed subgroups split any $n\mathbb{Z}$.
So, ${\mathfrak s}_{I}$ is undefined. On the other hand it is worth noting
that for every $n \in \omega$, $n \mathbb{Z} =_{a} Sym(\omega)$.
Indeed, fixing some representatives $a_{i}$ of all
the orbits, add the transpositions of the pairs $a_{i}, a_{i} + n$.
This induces all permutations on every orbit. 
Adding transpositions of some pairs from distinct orbits we get $Sym(\omega)$.  \parskip0pt

Another easy observation is that ${\mathfrak r}_{I} = 1$ in this case.
Indeed, the Pr\"{u}fer group $\mathbb{Z}(p^{\infty})$ with the natural action
on itself forms a reaping family.

\bigskip

It is interesting to compare the lattices that we consider here with the
lattice $P(\omega)$ of all subsets of $\omega$ and the ideal of finite
subsets. Since $=_{a}$ is a congruence of $P(\omega)$,
the orthogonality of infinite $a$ and $b$ means the absence of $c$
such that $c \le_{a} a$ and $c \le_{a} b$. So the van Douwen's cardinals
can be defined only in terms of $\le_{a}$ (and originally it was so).
On the other hand, this does not hold in lattices of subgroups of
$Sym(\omega)$. 
Indeed, let $\sigma$ be a transposition of some pair in
$\omega$. 
Then $\mathbb{Z}^{\sigma}$ induces a closed subgroup of
$Sym(\omega)$ which is $a$-equivalent to $\mathbb{Z}$ with the above
action. 
Clearly, the intersection of these groups is trivial. \parskip0pt

In the case of $(LF,IF)$ the corresponding example is as follows.
Let infinite $A,B,C \subset \omega$ define a partition of $\omega$
and $R$ be a bijection between $A$ and $B$. Let
$E_{0} = A^{2} \cup B^{2} \cup id_{C \times C}$ and
$E_{1} = R \cup id_{C \times C}$. It is easily seen that $E_{0}$
and $E_{1}$ are orthogonal equivalence relations, but $E_{1} \le_{a} E_{0}$.
The groups $G_{E_{0}} \cap SF(\omega)$ and $G_{E_{1}} \cap SF(\omega)$
have the same properties.

\end{document}